\documentclass[12pt]{amsart}

\usepackage[english]{babel}

\usepackage{amsmath}

\usepackage{changes}

\usepackage{amsthm}

\usepackage{amsfonts} 

\usepackage{hyperref}

\usepackage{epsfig}

\usepackage{amssymb}

\usepackage{environ}

\usepackage{mathrsfs}

\usepackage{caption}

\usepackage{tikzscale}

\usepackage{tikz-cd}

\usetikzlibrary{cd}

\usetikzlibrary{patterns}

\usepackage{enumitem}

\usepackage{dsfont}

\usepackage{aligned-overset}

\usepackage{mathtools}

\usepackage{float}

\usepackage{subcaption}

\usepackage{yfonts}

\usepackage[
backend=biber,
style=alphabetic,
]{biblatex}

\usepackage{import}

\tikzcdset{scale cd/.style={every label/.append style={scale=#1},
    cells={nodes={scale=#1}}}}

\DeclareMathOperator{\krn}{Ker}

\newtheorem{theorem}{Theorem}[section]
\newtheorem{definition}[theorem]{Definition}

\newtheorem{proposition}[theorem]{Proposition}

\newtheorem{remark}[theorem]{Remark}

\theoremstyle{definition}

\newcommand{\Ocal}{\mathcal{O}}

\DeclareMathAlphabet{\mathcal}{OMS}{cmsy}{m}{n}
\DeclareMathOperator{\tor}{Tor}
\usepackage{fullpage}

\title{A global proof of the homological excess intersection formula}
\author{Oscar Finegan}
\begin{document}

\maketitle
\begin{abstract}
    We provide a novel proof of the homological excess intersection formula for local complete intersections. The novelty is that the proof makes use of global morphisms comparing the intersections to a self intersection.
\end{abstract}
\section{Introduction}

The purpose of this short paper is to provide a novel proof of an existing result which can be found in \cite{https://doi.org/10.48550/arxiv.1510.04889}, namely, the homological excess intersection formula for local complete intersections:
\begin{theorem}\label{Theorem: Intro}
    Let $X$ be a nonsingular variety over an algebraically closed field $k$ of characteristic 0. Let $Y_1,\dots,Y_n$ be local complete intersection subvarieties such that their intersection $W$ is a local complete intersection in $X$. Then 
    \[Tor_q^{\Ocal_X}(\Ocal_{Y_1},\dots,\Ocal_{Y_n}) \cong \bigwedge^q\mathcal{E}_W\]
    where $\mathcal{E}_W = \krn\left(\bigoplus(\mathcal{N}_{Y_i/X}^\vee)|_W \rightarrow \mathcal{N}_{W/X}^\vee\right)$ is the excess bundle on $W$, which is of rank equal to the excess codimension $(\sum_{i = 1}^n codim(Y_i,X)) - codim(W,X).$ 
\end{theorem}

We are ultimately interested in extending this result to the case where the intersection of the $Y_i$ is no longer a local complete intersection. The strategy of the original proof is roughly to construct local isomorphisms using algebraic computations with Koszul complexes, and then to prove that the local isomorphisms glue up to a global one. When the intersection of the $Y_i$ is not itself lci, one runs into difficulties with computing Koszul cohomologies as well as with the gluing argument. In this work, we give a global morphism from the derived intersection of the $Y_i$ in $X$ to the derived self-intersection of the product of the $Y_i$ inside the $n$-fold product of $X$. Here, by derived intersection we mean the object $\Ocal_{Y_1}\otimes^L \dots \otimes^L \Ocal_{Y_n} \in \textbf{D}(X)$ etc.. The morphism is the unit for the pullback-pushforward adjunction of the closed immersion of the intersection in $X$. This morphism fits the multitors into a long exact sequence in which we have an explicit description for each third term \textit{with no assumptions on the intersection at all}. In the situation of Theorem \ref{Theorem: Intro}, this long exact sequence simplifies enough to prove the result.

\subsection*{Notation}
In this paper, $k$ will be an algebraically closed field of characteristic 0. 
\section{Preliminaries}

\subsection{Koszul complexes and local complete intersections}
Much of the material in this section is from \cite{F-L} and \cite{Matsumura}.
Let $X$ be a scheme, $\mathcal{E}$ a locally free $\Ocal_X$-module and $s : \mathcal{E} \rightarrow \Ocal_X$ a map of $\Ocal_X$-modules. 
\begin{definition}
    The Koszul complex $K_\bullet(\mathcal{E},s)$ is a differential graded algebra. As a graded algebra it is given by $\bigwedge^\bullet \mathcal{E}$ and the differential $d: \bigwedge^\bullet\mathcal{E} \rightarrow \bigwedge^\bullet\mathcal{E}$ is given as the unique derivation such that on degree one elements $d(e) = s(e)$. Locally on sections it is given by 
    \[e_1\wedge \dots \wedge e_p \mapsto \sum_{i=1}^p (-1)^is(e_i)e_1\wedge \dots \wedge\hat{e_i}\wedge \dots \wedge e_p,\]
    where $\hat{e_i}$ means `skip' $e_i$.
\end{definition}
\begin{remark}
    The uniqueness of the derivation here comes from the universal property for the tensor algebra, and the fact that this explicit local description of the differential kills elements of the form $e\wedge e$.
\end{remark}
There is an analogous algebraic definition of Koszul complex, and geometric Koszul complexes are locally isomorphic to algebraic Koszul complexes. In the algebraic case, if one chooses a basis for $\mathcal{E}$ and $f_i$ is the image of the $i^{th}$ basis element under $s$, we write $K_\bullet(f_1,\dots,f_n)$ for $K_\bullet(\mathcal{E},s)$. Note that picking different basis elements of $\mathcal{E}$ will result in isomorphic Koszul complexes, with the isomorphism given by a change of basis. We are interested in these complexes for two reasons. Firstly, they behave nicely with respect to the tensor product of complexes.

\begin{proposition}
Suppose we have two morphisms $s: \mathcal{E} \rightarrow \Ocal_X, t : \mathcal{F} \rightarrow \Ocal_X$ of locally free $\Ocal_X$-modules of finite rank. Then we have an identification 
\[K^\bullet(\mathcal{E},s)\otimes K^\bullet(\mathcal{F},t) \cong K^\bullet(\mathcal{E}\oplus \mathcal{F},s\oplus t).\]
\end{proposition}

The second reason we are interested in Koszul complexes is that they locally give free resolutions for structure sheaves of local complete intersections. 
To demonstrate this we first need to introduce some notions of regularity. Let $R$ be a ring (always assumed commutative and unital). We say that a sequence $f_1,\dots,f_n \in R$ is regular if $f_i$ is not a zero-divisor in the quotient ring $R/(f_1,\dots,f_{i-1})$. We say that an ideal sheaf $\mathcal{I} \subset \Ocal_X$ is a regular ideal sheaf if it can locally be generated by a regular sequence of sections of $\Ocal_X$. We call a section $s: \mathcal{E} \rightarrow \Ocal_X$ a regular section if its image sheaf is a regular ideal sheaf. Any regular ideal sheaf is locally the image sheaf of a regular section. Indeed, working affine locally, the ideal $(f_1,\dots,f_n)$ is the image of the map $R^n \rightarrow R$ sending $e_i \mapsto f_i$. We have the following relation between Koszul complexes and regular sections:

\begin{proposition}\cite{F-L}
For any regular section $s: \mathcal{E} \rightarrow \Ocal_X$, $K^\bullet(\mathcal{E},s)$ is a global locally free resolution of the structure sheaf of $Z(s)$, the zero locus of the section $s$.
\end{proposition}
In the case that the scheme $X$ is noetherian, this proposition is reversible, i.e. the vanishing of Koszul cohomologies implies that the section $s$ is regular \cite{Matsumura}.
A closed subscheme $Y \subset X$ is called a \textit{local complete intersection} if locally around every closed point $y \in Y$, there is a neighbourhood $U$ of $y$ in $X$ such that $\mathcal{I}_Y(U)$ is generated by codim($Y,X$) elements. In the case where $X$ is a nonsingular variety over $k$, we make use of the following 

\begin{theorem}[\cite{Matsumura} Theorem 17.4]\label{Matsumura}
    Let $(A,\mathfrak{m})$ be a noetherian Cohen-Macaulay local ring. Then 
    \begin{enumerate}
        \item For every ideal $I$ in $A$, we have an equality $ht(I) + dim(A/I) = dim(A)$.
        \item  A sequence $a_1,\dots,a_r \in \mathfrak{m}$ is regular $\iff ht(a_1,\dots,a_r) = r.$
    \end{enumerate}
\end{theorem}

This shows that a local complete intersection is defined by a regular ideal sheaf in a noetherian Cohen-Macaulay scheme, and therefore locally has Koszul resolutions of its structure sheaf. 
A local complete intersection can be seen as the generalisation of the notion of effecive Cartier divisor to higher codimensions. Indeed, effective Cartier divisors are defined by having regular ideal sheaves locally generated by a single regular element. An important example is supplied by the fact that any nonsingular subvariety of a nonsingular variety is a local complete intersection.

\subsection{Tor-independence}
The results and definitions in this section are from \cite{Lipman}. All of the functors in this section are derived and we suppress the $R,L$ that usually attends derived functors. In the proofs of the main results later, we have fibre squares 
\[\begin{tikzcd}
        &X' \ar[r,"v"]\ar[d,"g"]\ar[dr,phantom,"\sigma"] & X\ar[d,"f"]\\
        &Y' \ar[r,"u"] & Y
    \end{tikzcd}\]
and we make use of the base change map $f^*u_* \rightarrow v_*g^*$ in certain cases in which it is an isomorphism. Hence we include a brief overview on the base-change map and the notion of Tor-independence. If we have any commuting square $\sigma$ as above (not necessarily fibre), we define the base change map 
\[\beta_\sigma : f^*u_* \rightarrow v_*g^*\]
to be the composition
\[f^*u_* \rightarrow f^*u_*g_*g^* \rightarrow f^*f_*v_*g^* \rightarrow v_*g^*\]
where the first and third arrows come from the pullback-pushforward adjunctions and the middle arrow from the commutativity of the square.
\begin{definition}\label{def:tor-ind}
    A fibre square of schemes over a scheme $S$
    \[\begin{tikzcd}
        &X' \ar[r,"v"]\ar[d,"g"] & X\ar[d,"f"]\\
        &Y' \ar[r,"u"] & Y
    \end{tikzcd}\]
    is said to be \textit{Tor-independent} if for all pairs of points $y' \in Y', x \in X$ such that $ u(y') = y = f(x):$
    \[Tor_q^{\Ocal_{Y,y}}(u_*\Ocal_{Y',y'},f_*\Ocal_{X,x}) = 0 \quad \forall q>0.\]
\end{definition}

One can see almost immediately from the definition that we have the following result on Tor-independent fibre squares.

\begin{proposition}\label{Prop: Flat}
    If the fibre square $\sigma$ in Definition \ref{def:tor-ind} has either $f$ or $u$ being a flat morphism, then it is Tor-independent. 
\end{proposition}

The benefit of Tor-independence is that it gives a homological criterion for the base change map to be an isomorphism in cases where the schemes and maps involved are relatively nice. 

\begin{definition}\label{def:concentrated}
A map of schemes $f: X \rightarrow Y$ is said to be concentrated if it is quasi-compact and quasi-separated.  
\end{definition}

\begin{theorem}[\cite{Lipman} Theorem 3.10.3]\label{prop:base change}
    For any fibre square as in Definition \ref{def:tor-ind} where the maps are concentrated and the schemes quasi-separated, Tor-independence is equivalent to the base change map for the square being an isomorphism;
    \[f^*u_* \cong v_*g^*\]
\end{theorem}

Since the results of this paper are about closed immersions of varieties over a field $k$, we will be able to make use of this theorem in the proofs of the main results. Namely, once we know a fibre square is Tor-independent, we may use that the base change map around that square is an isomorphism. 

We need to use the base change isomorphism later for an equivalence computation, and here we include a proof of a result saying that a fibre square for the proper intersection of local complete intersections whose intersection is also a local complete intersection is Tor-independent. 
\begin{proposition}\label{Prop: Proper intersection is tor independent}
    Let $X$ be a Cohen-Macaulay scheme with local complete intersection subschemes $Y$ and $Z$. If $Y \cap Z$ has the expected codimension $(codim(Y\cap Z) = codim(Y) + codim(Z))$ then $Y \cap Z$ is a local complete intersection in $X$ and the following fibre square of closed immersions is Tor-independent
    \[\begin{tikzcd}
        &Y\cap Z \ar[r]\ar[d] &Y \ar[d]\\
        &Z \ar[r] &X
    \end{tikzcd}\]
\end{proposition}

\begin{proof}
The vanishing of the $Tor_q$ is a local criterion so let us work affine locally. Let the local defining ideals for $Y$ and $Z$ be given by $(f_1,\dots,f_n)$ and $(g_1,\dots,g_m)$ respectively. The ideal $I = (f_1,\dots,f_n,g_1,\dots,g_m)$ defines $Y\cap Z$ so $Y\cap Z$ is an lci. From Theorem \ref{Matsumura}, we know that $ht(I) = n+m$ and also that the sequence $f_1,\dots, f_n,g_1,\dots, g_m$ is regular. The $Tor_q^X(Y,Z)$ can be computed locally by tensoring the Koszul resolutions for the structure sheaves, i.e. locally it is given by the cohomology of $K^\bullet(f_1,\dots,f_n,g_1,\dots,g_m)$. These cohomologies vanish for $q>1$ since we have just shown that this sequence is regular. 
\end{proof}

This Proposition can be interpreted as saying that proper intersections of local complete intersections have no excess.

\section{Results}

In this section we provide proofs for results that exist in the literature on the geometric properties of the derived intersections of local complete intersections whose intersection is also a local complete intersection. We will first give a classical account of the derived self-intersection formula for a local complete intersection. We will need to use this result and the construction of the isomorphism in the proof of our main result. The second part of this section is dedicated to a novel proof of a result by Scala in \cite{https://doi.org/10.48550/arxiv.1510.04889} of an excess intersection formula for lci intersections.

To begin, we prove that within the framework of intersections of local complete intersections the local computations of Koszul cohomologies will glue to agree with a global model for the multitors. Let $Y_1,\dots, Y_n$ be local complete intersection subvarieties of a scheme $X$. Pick global flat resolutions $\mathcal{F}_i^\bullet \rightarrow \Ocal_{Y_i}$ for each $i$. We fix a global model of the multitors by defining 
\[\tor_q^{\Ocal_X}(\Ocal_{Y_1},\dots, \Ocal_{Y_n}) := H^{-q}(\mathcal{F}^\bullet_1 \otimes \dots \otimes \mathcal{F}^\bullet_n).\]
For $p \in Z := \bigcap Y_i$, let $U$ be a local affine neighbourhood such that each $\Ocal_{Y_i}|_U$ can be resolved by a free Koszul complex $K^\bullet(E_i,s_i)$, which exists because each $Y_i$ is a local complete intersection. Restriction to $U$ is an exact functor which preserves flatness, so the restriction of the augmentation $\mathcal{F}^\bullet_i|_U \rightarrow \Ocal_{Y_i}|_U$ is still a quasi-isomorphism. Then since the complexes $K^\bullet(E_i,s_i)$ are free resolutions of $\Ocal_{Y_i}|_U$ (and therefore q-projective), there are morphisms (unique up to homotopy) $\psi_i: K^\bullet(E_i,s_i) \rightarrow \mathcal{F}_i^\bullet|_U$ lifting the augmentation to $\Ocal_{Y_i}|_U$. That is, making the following triangle commute;
\[
\begin{tikzcd}[every arrow/.append style={shift left}]
    &&\mathcal{F}_i^\bullet|_U \ar[d] \\
    &K^\bullet(E_i,s_i) \ar[r]\ar[ur, dashed, "\psi_i"] &\Ocal_{Y_i}|_U
\end{tikzcd}
\]
Necessarily, the $\psi_i$ are quasi-isomorphisms because both of the map from $K^\bullet(E_i,s_i)$ and $\mathcal{F}_i^\bullet|_U$ to $\Ocal_{Y_i}|_U$ are quasi-isomorphisms. These lifts $\psi_i$ together induce a quasi-isomorphism 
\[\psi:  K^\bullet(E_1,s_1) \otimes \dots \otimes K^\bullet(E_n,s_n) \cong K^\bullet\left(\bigoplus E_j, \oplus s_j\right) \rightarrow \mathcal{F}_1^\bullet|_U \otimes \dots \otimes \mathcal{F}_n^\bullet|_U,\]
because each $\mathcal{F}_i^\bullet|_U$ is a q-flat complex and tensoring with q-flat complexes preserves quasi-isomorphisms. This quasi-isomorphism $\psi$ is the unique morphism (up to homotopy) lifting the natural map $\bigotimes \mathcal{F}_i^\bullet|_U \rightarrow \bigotimes \Ocal_{Y_i}|_U$. Since $\psi$ is a quasi-isomorphism, it gives isomorphisms in cohomology 
\[H^{-q}(\psi): H^{-q}\left(K^\bullet\left(\bigoplus E_j, \oplus s_j\right)\right) \rightarrow H^{-q}(\mathcal{F}_1^\bullet|_U \otimes \dots \otimes \mathcal{F}_n^\bullet|_U) \cong Tor_q^{\Ocal_X}(\Ocal_{Y_1},\dots, \Ocal_{Y_n})|_U .\]
As the $\psi_i$ are unique up to homotopy, the $H^{-q}\psi$ are unique isomorphisms on cohomology coming from lifts of the augmentations to $\Ocal_{Y_i}|_U$.
\par
Take another open set $U' \subset X$ such that $U \cap U' = V \neq \emptyset$ and such that there are free Koszul resolutions $K^\bullet(E_i',s_i') \rightarrow \Ocal_{Y_i}|_{U'}$.
On $V$, we have a commuting diagram of lifts unique up to homotopy for each $\Ocal_{Y_i}|_V$
\[\begin{tikzcd}
    &K^\bullet(E_i,s_i)|_V \ar[dr] \ar[r,"\psi_i|_V"]&\mathcal{F}^\bullet_i|_V\ar[d]&K^\bullet(E'_i,s'_i)|_V\ar[dl]\ar[l,"\psi_i'|_V"']\\
    &&\Ocal_{Y_i}|_V.&
\end{tikzcd}\]
There are therefore diagrams commuting up to homotopy; 
\[\begin{tikzcd}[column sep = 0.5em]
    &K^\bullet(E_i,s_i)|_V \ar[rr,"\chi_i"] \ar[dr,"\psi_i|_V"'] &&K^\bullet(E_i',s_i')|_V \ar[dl,"\psi_i'|_V"]\\
    &&\mathcal{F}^\bullet_i|_V &\\
\end{tikzcd}\]
of lifts of the augmentations to $\Ocal_{Y_i}|_V$. Hence, setting $E:= \bigoplus E_i, E':= \bigoplus E_i'$, the induced quasi-isomorphisms
\[\begin{tikzcd}[column sep = 0.5em]
    &K^\bullet(E,s)|_V \ar[rr,"\chi"] \ar[dr,"\psi|_V"'] &&K^\bullet(E',s')|_V \ar[dl,"\psi'|_V"]\\
    &&\bigotimes \mathcal{F}_i^\bullet|_V &\\
\end{tikzcd}\]
also commute up to homotopy and so the induced isomorphisms 
\[\begin{tikzcd}[column sep = 0.5em]
    &H^{-q}(K^\bullet(E,s)|_V) \ar[rr,"H^{-q}\chi"] \ar[dr,"H^{-q}\psi |_V"'] &&H^{-q}(K^\bullet(E',s')|_V) \ar[dl,"H^{-q}\psi' |_V"]\\
    &&\tor_q(\Ocal_{Y_1},\dots,\Ocal_{Y_n})|_V &\\
\end{tikzcd}\]
commute on the nose. The commutativity of these triangles is what we mean when we talk about the fact that the cohomologies of the local Koszul models glue.  
\par
For us there are natural choices of representative for the homotopy equivalences $\chi_i : K^\bullet(E_i,s_i)|_V \rightarrow K^\bullet(E_i',s_i')|_V$. Namely, since the ideal generated by the images of $s_i|_V$ and $s_i'|_V$ are the same, there is a change of basis morphism $\chi_{i,0} : E_i|_V \rightarrow E_i'|_V$ making the triangle 
\[\begin{tikzcd}[column sep = 0.5em]
    &E_i|_V \ar[rr,"\chi_{i,0}"] \ar[dr,"s_i|_V"'] &&E_i'|_V \ar[dl,"s_i'|_V"]\\
    &&\Ocal_V &\\
\end{tikzcd}\]
commute. This change of basis morphism extends by the functoriality of the Koszul complex to a morphism of Koszul complexes lifting the augmentations to $\Ocal_{Y_i}|_V$, and hence by the uniqueness property from q-projectivity is a homotopy equivalence. 

\begin{theorem}\label{Theorem: Self-intersection}
    Let $i: Y \rightarrow X$ be a local complete intersection subscheme of a nonsingular variety over $k$. Then
    \[\tor_q(i_*\Ocal_Y,i_*\Ocal_Y) \cong i_*\bigwedge^q \mathcal{C}_{Y/X}\]
    where $\mathcal{C}_{Y/X} = \mathcal{N}_{Y/X}^\vee$ is the conormal bundle for $Y$ in $X$.
\end{theorem}
\begin{proof}
Let $\{U_\alpha\}$ be an open cover of $X$ such that on each $U_\alpha$ there is a Koszul resolution $K^\bullet(E_\alpha,s_\alpha)$ of $i_*\Ocal_Y|_{U_\alpha}$. For each $\alpha$ we have the following Tor-independent fibre square of closed and open immersions (Proposition \ref{Prop: Flat})
\[
\begin{tikzcd}
    &Y \cap U_\alpha \ar[r,"i_\alpha"]\ar[d,"(j_\alpha)_Y"']& U_\alpha \ar[d,"j_\alpha"]\\
    &Y \ar[r,"i"]& X.
\end{tikzcd}
\]
 We want to construct local isomorphisms 
\[H^{-q}(K^\bullet(E_\alpha,s_\alpha)\otimes (i_*\Ocal_Y)|_{U_\alpha}) \rightarrow (i_*\bigwedge \mathcal{C}_{Y/X})|_{U_\alpha}\]
for each $\alpha$, commuting with the isomorphisms induced by the change of basis morphisms

\[
\begin{tikzcd}[column sep = 0.3em]
    &H^{-q}(K^\bullet(E_\alpha,s_\alpha)|_{U_{\alpha\beta}}\otimes (i_*\Ocal_Y)|_{U_{\alpha\beta}})\ar[dr]\ar[rr] &&H^{-q}(K^\bullet(E_\beta,s_\beta)|_{U_{\alpha\beta}}\otimes (i_*\Ocal_Y)|_{U_{\alpha\beta}})\ar[dl]\\
    &&(i_*\bigwedge \mathcal{C}_{Y/X})|_{U_{\alpha\beta}},& 
\end{tikzcd}
\]
where $U_{\alpha\beta}$ means $U_{\alpha} \cap U_\beta$. Doing so will prove the statement of the theorem, because the discussion above tells us that there are local isomorphisms of the Tors to Koszul models which commute with the isomorphism on cohomology coming from the change of basis map. For each $\alpha$, the complex $K^\bullet(E_\alpha, s_\alpha) \otimes (i_*\Ocal_Y)|_{U_\alpha}$ has zero differential  and has $H^{-q}(K^\bullet(E_\alpha,s_\alpha)\otimes (i_*\Ocal_Y)|_{U_\alpha}) = K^{-q}(E_\alpha,s_\alpha)\otimes (i_*\Ocal_Y)|_{U_\alpha} \cong (i_\alpha)_*\bigwedge^qE_\alpha|_{Y\cap U_\alpha} \cong (i_\alpha)_*\bigwedge^q\mathcal{C}_{Y\cap U_\alpha/U_\alpha}$. 
This final isomorphism is the restriction of the morphism $s_\alpha$ to $Y \cap U_\alpha$, which becomes an isomorphism onto its image. We now claim that 
\[(i_\alpha)_*\bigwedge^q \mathcal{C}_{Y \cap U_\alpha/U_\alpha} \cong (j_\alpha)^*i_*\bigwedge^q \mathcal{C}_{Y/X}.\]
To see this, note that $\mathcal{C}_{Y/X}:= i^*\mathcal{I}_Y$, $\mathcal{C}_{Y \cap U_\alpha/U_\alpha}:= (i_\alpha)^*\mathcal{I}_{Y\cap U_\alpha} = (j_\alpha \circ i_\alpha)^*\mathcal{I}_Y$. Then there are isomorphisms 
\[(i_\alpha)_*\bigwedge^q(j_\alpha \circ i_\alpha)^*\mathcal{I}_Y \cong (i_\alpha)_*\bigwedge^q(i \circ (j_\alpha)_Y)^*\mathcal{I}_Y 
\cong (i_\alpha)_*(j_\alpha)_Y^*\bigwedge^q\mathcal{C}_{Y/X}
\cong (j_\alpha)^*i_*\bigwedge^q\mathcal{C}_{Y/X},\]
where the first isomorphism comes from commutativity of the above square, the second from commutativity of exterior powers with pullbacks, and the third from base change around the square. Then on $U_{\alpha\beta}$ we have the following diagram of isomorphisms; 
\[
\begin{tikzcd}[column sep = 0.2em]
    &H^{-q}(K^\bullet(E_\alpha,s_\alpha)|_{U_{\alpha\beta}}\otimes (i_*\Ocal_Y)|_{U_{\alpha\beta}})\ar[d]\ar[rr] &&H^{-q}(K^\bullet(E_\beta,s_\beta)|_{U_{\alpha\beta}}\otimes (i_*\Ocal_Y)|_{U_{\alpha\beta}})\ar[d]\\
    &(i_{\alpha\beta})_*(\bigwedge^q\mathcal{C}_{Y\cap U_\alpha/U_\alpha})|_{U_{\alpha\beta}} \ar[dr]&&(i_{\alpha\beta})_*(\bigwedge^q\mathcal{C}_{Y\cap U_\beta/U_\beta})|_{U_{\alpha\beta}} \ar[dl]\\
    &&(i_*\bigwedge^q \mathcal{C}_{Y/X})|_{U_{\alpha\beta}}.&
\end{tikzcd}
\]
 Since the diagrams
\[
\begin{tikzcd}[column sep = 0.2em]
    &E_\alpha|_{U_{\alpha\beta}}\ar[dr,"s_\alpha|_{U_{\alpha\beta}}"']\ar[rr]&&E_\beta|_{U_{\alpha\beta}}\ar[dl,"s_\beta|_{U_{\alpha\beta}}"] \\
    &&\Ocal_{U_{\alpha\beta}},&
\end{tikzcd}
\]
commute, where the top morphism is the change of basis morphism,
and the vertical isomorphisms in the pentagonal diagram come from the restrictions of the $s_\alpha$, we see the diagram commutes. This is exactly what we were trying to show and we conclude that the local isomorphisms glue to give a global isomorphism
\[\tor^{\Ocal_X}_q(i_*\Ocal_Y,i_*\Ocal_Y) \cong i_*\bigwedge^q\mathcal{C}_{Y/X}.\]
\end{proof}

We need to change our perspective slightly on our multitors in order for our new approach to work.

\begin{proposition}\label{prop: equivalence}
    Let $X$ be a non-singular scheme over $k$ with local complete intersection subschemes $i_j: Y_j \rightarrow X$ $(1\leq j\leq n)$, whose intersection we denote by $W$. Consider the fiber square 
    \[\begin{tikzcd}[column sep = 2.5em]
        &W \arrow[r,"w"]\arrow[d,"j_1\times \dots \times j_n"']& X \arrow[d,"\Delta_X^n"]\\
        &Y_1\times \dots \times Y_n \arrow[r,"(i_1\times \dots \times i_n)"]& X^{\times n}
    \end{tikzcd}\]
    where $X^{\times n}$ means the $n$-fold fibre product of $X$ with itself over $k$, $\Delta_X^n$ is the diagonal morphism, and $w$ is the closed embedding of $W$ into $X$. Suppressing the notation of right and left derived functors, there is an isomorphism in $\textbf{D}(X)$;
    \[(i_1)_*\Ocal_{Y_1}\otimes \dots \otimes (i_n)_*\Ocal_{Y_n} \cong (\Delta_X^n)^*(i_1\times \dots \times i_n)_*\Ocal_{Y_1\times \dots \times Y_n}.\]
\end{proposition}

\begin{proof}
     We prove the statement by induction on $n$. The case $n = 1$ is clear so  
    suppose $n > 1$.  Consider the expanded commutative diagram;
\begin{center}
\begin{tikzcd}[column sep = 5em]
        &W \arrow[r,"j_n"]\arrow[dd,"j_1\times \dots \times j_n"']&Y_n \ar[r, "i_n"]\ar[d, "\Delta^n_{Y_n}"]  &X \arrow[dd,"\Delta_X^n"]\\
        &&Y_n^{\times n}\ar[d,"(i_n\times\dots \times i_n)\times 1"] &\\
        &(Y_1\times \dots \times Y_{n-1})\times Y_n \arrow[r,"(i_1\times \dots \times i_{n-1}) \times 1"]& X^{\times n-1}\times Y_n \ar[r,"(1\times \dots \times 1)\times i_n"] &X^{\times n}
    \end{tikzcd}
\end{center}
    where $\Delta^n_{Y_n}$ is the diagonal morphism for $Y_n$ into $Y_n^{\times n}$. Here $\Delta_X^n$ is a regular immersion because $X$ is non-singular. By Proposition \ref{Prop: Proper intersection is tor independent} the square on the right hand side is Tor-independent. So then
    \[(\Delta^n_X)^*(i_1\times \dots \times i_n)_* = (i_n)_*(\Delta^n_{Y_n})^*(i_n\times \dots \times  i_n\times 1)^*(i_1\times\dots \times i_{n-1}\times 1)_*.\]
    Denote by $\pi_{(1,\dots,n-1)}$ the projection to the first $n-1$ factors. Then since $\Ocal_{Y_1\times \dots \times Y_n} = \pi_{(1,\dots,n-1)}^*\Ocal_{Y_1\times \dots \times Y_{n-1}}$ we have the identification
    \begin{align*}
        &(i_n)_*(\Delta^n_{Y_n})^*(i_n\times \dots \times  i_n\times 1)^*(i_1\times\dots \times i_{n-1}\times 1)_*(\Ocal_{Y_1\times \dots \times Y_n})\\
        &=(i_n)_*(\Delta^n_{Y_n})^*(i_n\times \dots \times  i_n\times 1)^*(i_1\times\dots \times i_{n-1}\times 1)_*\pi_{(1,\dots,n-1)}^*(\Ocal_{Y_1\times \dots \times Y_{n-1}}).
    \end{align*}
    There is another commuting diagram;
    \begin{center}
        \begin{tikzcd}[column sep = 5em]
        &(Y_1\times \dots \times Y_{n-1})\times Y_n \ar[r,"(i_1\times \dots \times i_{n-1})\times 1"]\ar[d, "\pi_{(1, \dots, n-1)}"] &X^{\times (n-1)}\times Y_n \ar[d,"\pi_{(1, \dots, n-1)}"] & Y_n^{\times n} \ar[l,"(i_n\times \dots \times i_n)\times 1"']\ar[d,"\pi_{(1,\dots,n-1)}"]\\
        &Y_1\times \dots \times Y_{n-1} \ar[r,"i_1\times \dots \times i_{n-1}"] &X^{\times (n-1)} &Y_n^{\times (n-1)}\ar[l,"i_n\times \dots \times i_n"']\\
        &\bigcap_{i = 1}^{n-1}Y_i \ar[u]\ar[r] &X \ar[u,"\Delta^{n-1}_X"'] & Y_n .\ar[l,"i_n"']\ar[u,"\Delta^{n-1}_{Y_n}"']
    \end{tikzcd}
    \end{center}
    Here now by Proposition \ref{Prop: Flat} and because projections are flat morphisms the top left hand square is Tor-independent. Hence 
    \[(i_n\times \dots \times  i_n\times 1)^*(i_1\times\dots \times i_{n-1}\times 1)_*\pi_{(1,\dots,n-1)}^* =  \pi_{(1,\dots,n-1)}^*(i_n\times \dots \times i_n)^*(i_1\times \dots \times i_{n-1})_*.\]
    Furthermore the composed morphism
    \[\begin{tikzcd}[column sep = 5em]
        &Y_n \ar[r,"\Delta_{Y_n}^n"] & Y_n^{\times n}\ar[r,"\pi_{(1,\dots,n-1)}"] &Y_n^{\times (n-1)}
    \end{tikzcd}\]
    and the morphism $\Delta_{Y_n}^{n-1}$ are equal. Hence 
    \[(\Delta^n_{Y_n})^*\pi_{(1,\dots,n-1)}^*(i_n\times \dots \times i_n)^*(i_1\times \dots \times i_{n-1})_* =(\Delta^{n-1}_{Y_n})^*(i_n\times \dots \times i_n)^*(i_1\times \dots \times i_{n-1})_*\]
    By commutativity of the bottom right square, we have an equality
    \[(\Delta^{n-1}_{Y_n})^*(i_n\times \dots \times i_n)^*(i_1\times \dots \times i_{n-1})_* = (i_n)^*(\Delta_X^{n-1})^*(i_1\times \dots \times i_{n-1})_*\]
    Combining all of the above identifications, we have shown that there is an isomorphism in $\textbf{D}(X)$
    \[(\Delta_X^n)^*(i_1\times \dots \times i_n)_*\Ocal_{Y_1\times \dots \times Y_n} \cong (i_n)_*(i_n)^*(\Delta_X^{n-1})^*(i_1\times \dots \times i_{n-1})_*(\Ocal_{Y_1\times \dots \times Y_{n-1}}).\]
    The result then follows from the projection formula and induction.   
\end{proof}

\begin{theorem}\label{Theorem: Excess Intersection Formula}
    Let $X$ be a non-singular variety over an algebraically closed field $k$ and let $Y_1,\dots,Y_n$ be local complete intersection subvarieties of $X$. Suppose that $w:W = \bigcap Y_i\rightarrow X$ is also a local complete intersection. Then
    \[Tor_q^{\Ocal_X}(\Ocal_{Y_1},\dots,\Ocal_{Y_n}) \cong w_*\bigwedge^q\mathcal{E}_W,\]
    where $\mathcal{E}_W$ (the excess bundle) is defined as the kernel of the natural surjection 
    \[\bigoplus \mathcal{N}^\vee_{Y_i/X}|_W \rightarrow \mathcal{N}_{W/X}^\vee\]
    coming from $\mathcal{I}_W = \sum \mathcal{I}_{Y_i}.$
\end{theorem}

\begin{proof}
    We prove this statement when $n =2$ for notational convenience, with the general case proceeding by an equivalent argument. By Proposition \ref{prop: equivalence} we are looking to compute the cohomologies of the object $\Delta^*(i_1 \times i_2)_*\Ocal_{Y_1\times Y_2}$. The closed immersion $w$ gives rise to an adjunction unit morphism
    \[\eta_w: \Delta^*(i_1\times i_2)_*\Ocal_{Y_1\times Y_2} \rightarrow w_*w^*\Delta^*(i_1 \times i_2)_*\Ocal_{Y_1\times Y_2}.\] 
    We have an isomorphism
    \[w_*w^*\Delta^*(i_1 \times i_2)_*\Ocal_{Y_1\times Y_2} \cong w_*(j_1\times j_2)^*(i_1 \times i_2)^*(i_1 \times i_2)_*\Ocal_{Y_1 \times Y_2},\]
    from commutativity of the square in the statement of Proposition \ref{prop: equivalence}. 
    We can compute the cohomologies of the object on the right, since the morphism $(i_1 \times i_2)$ is a regular immersion. Applying the derived self-intersection formula Theorem \ref{Theorem: Self-intersection}, and writing $\mathcal{C}_{Y/X}$ for $\mathcal{N}_{Y/X}^\vee$ etc,we see that there are global isomorphisms
    \[H^{-q}((i_1\times i_2)^*(i_1 \times i_2)_*\Ocal_{Y_1\times Y_2}) \cong \bigwedge^q\mathcal{C}_{Y_1\times Y_2/X\times X}.\]
    Since these cohomologies are locally free sheaves on $Y_1\times Y_2$, taking cohomology commutes with the pullback $(j_1\times j_2)^*$. Additionally $w$ is a closed immersion so $w_*$ is an exact functor and also commutes with cohomology. We therefore compute that there is a global isomorphism
    \[\varphi_q: H^{-q}(w_*w^*\Delta^*(i_1 \times i_2)_*\Ocal_{Y_1\times Y_2}) \rightarrow w_*(j_1\times j_2)^*\bigwedge^q\mathcal{C}_{Y_1\times Y_2/X\times X} \cong w_*\bigwedge^q(\mathcal{C}_{Y_1/X}\oplus \mathcal{C}_{Y_2/X})|_W. \]
    This isomorphism has the following important property. If $U$ is any open subset on which $\Ocal_{Y_1}$ and $\Ocal_{Y_2}$ have Koszul resolutions $K^\bullet(\mathcal{F}_1, s_1), K^\bullet(\mathcal{F}_2,s_2)$ respectively, then $\varphi_q|_U$ factors as 
    \[H^{-q}(w_*w^*\Delta^*(i_1\times i_2)_*\Ocal_{Y_1\times Y_2})|_U \rightarrow w_*H^{-q}(K^\bullet(\mathcal{F}_1\oplus \mathcal{F}_2,s_1\oplus s_2)|_{W\cap U}) \rightarrow w_*\bigwedge^q(\mathcal{C}_{Y_1/X}\oplus \mathcal{C}_{Y_2/X})|_{W\cap U}, \]
    where the first morphism is the inverse of the one induced from the $K^\bullet(\mathcal{F}_i,s_i)$ being free resolutions and the second morphism coming from the computation of Koszul cohomology in this case.
    Since $\mathcal{E}_W$ is a subbundle of $\mathcal{C}_{Y_1/X}\oplus \mathcal{C}_{Y_2/X}$ we have a diagram of global morphisms 
    \[\begin{tikzcd}
        &H^{-q}(\Delta^*(i_1\times i_2)_*\Ocal_{Y_1\times Y_2}) \ar[r,"H^{-q}\eta_w"] & H^{-q}(w_*w^*\Delta^*(i_1\times i_2)_*\Ocal_{Y_1\times Y_2}) \ar[d,"\varphi_q"]\\
        &w_*\bigwedge^q\mathcal{E}_W \ar[r,"\iota"]& w_*\bigwedge^q(\mathcal{C}_{Y_1/X}\oplus \mathcal{C}_{Y_2/X})|_W.
    \end{tikzcd}\]
    with the bottom arrow being the natural inclusion and $\varphi_q$ the isomorphism above. Now all we need to check is that the image of $w_*\bigwedge^q\mathcal{E}_W$ under $\varphi_q^{-1}\circ \iota$ agrees with the image of $H^{-q}(\Delta^*(i_1\times i_2)_*\Ocal_{Y_1\times Y_2})$ under $H^{-q}\eta_w$ and that $H^{-q}\eta_w$ is an injection to conclude the statement of the theorem. Checking that two subsheaves are equal and that a map is an injection are both local checks so we may work locally. 
    \par
    We work locally enough that $\Ocal_{Y_1}$ has resolution $K^\bullet(\mathcal{F}_1,s_1)$, $\Ocal_{Y_2}$ has resolution $K^\bullet(\mathcal{F}_2,s_2)$ and $\Ocal_W$ has resolution $K^\bullet(\mathcal{G},t)$. Then, since the intersection of $Y_1$ and $Y_2$ is equal to $W$, we have 
    \[K^\bullet(\mathcal{F}_1,s_1)\otimes K^\bullet(\mathcal{F}_2,s_2) \cong K^\bullet(\mathcal{F}_1\oplus \mathcal{F}_2, s_1\oplus s_2) \cong K^\bullet(\mathcal{G}\oplus \mathcal{F}, t \oplus 0).\]
    We can compute that $H^{-q}(K^\bullet(\mathcal{G}\oplus \mathcal{F}, t \oplus 0) ) \cong \bigwedge^q\mathcal{F}|_W$. We assume that we are working locally enough that the short exact sequence of locally free sheaves on $W$
    \[0 \rightarrow \mathcal{E}_W \rightarrow (\mathcal{C}_{Y_1/X}\oplus \mathcal{C}_{Y_2/X})|_W \rightarrow \mathcal{C}_{W/X}\rightarrow 0\]
    splits. Then we have two split exact sequences
    \[0 \rightarrow \mathcal{F}|_W \rightarrow (\mathcal{F}_1 \oplus \mathcal{F}_2)|_W \rightarrow \mathcal{G}|_W \rightarrow 0\]
    and
    \[0 \rightarrow \mathcal{E}_W \rightarrow (\mathcal{C}_{Y_1/X}\oplus \mathcal{C}_{Y_2/X})|_W \rightarrow \mathcal{C}_{W/X}\rightarrow 0\]
    which have isomorphic middle and rightmost terms. Therefore they are isomorphic split exact sequences. In particular, we have a commuting square
    \[
    \begin{tikzcd}
        &\mathcal{F}|_W \ar[r]\ar[d,"\cong"] &(\mathcal{F}_1\oplus \mathcal{F}_2)|_W\ar[d,"\cong"]\\
        &\mathcal{E}_W \ar[r]& (\mathcal{C}_{Y_1/X}\oplus \mathcal{C}_{Y_2/X})|_W.
    \end{tikzcd}
    \]
    We want now to show that $H^{-q}\eta_w$ is an injection. To do this, we recall that the functor $w_*w^* \cong - \otimes \Ocal_W$. Then we get a long exact sequence on cohomology
    \begin{center}
    \begin{tikzpicture}[descr/.style={fill=white,inner sep=1.5pt}]
        \matrix (m) [
            matrix of math nodes,
            row sep=1em,
            column sep=2.5em,
            text height=1.5ex, text depth=0.25ex
        ]
        { \dots & H^{-q}(\Delta^*(i_1\times i_2)_*\Ocal_{Y_1\times Y_2}\otimes \mathcal{I}_W) & H^{-q}(\Delta^*(i_1\times i_2)_*\Ocal_{Y_1\times Y_2}) \\
            & H^{-q}(w_*w^*\Delta^*(i_1\times i_2)_*\Ocal_{Y_1\times Y_2}) &H^{-q + 1}(\Delta^*(i_1\times i_2)_*\Ocal_{Y_1\times Y_2}\otimes \mathcal{I}_W) & \dots \\
        };
        \path[overlay,->, font=\scriptsize,>=latex]
        (m-1-1) edge (m-1-2)
        (m-1-2) edge (m-1-3)
        (m-1-3) edge[out=355,in=175] node[descr,yshift=0.3ex] {$H^{-q}\eta_w$} (m-2-2)
        (m-2-2) edge (m-2-3)
        (m-2-3) edge (m-2-4);
    \end{tikzpicture}
    \end{center}
    coming from the short exact sequence
    \[0 \rightarrow \mathcal{I}_W \rightarrow \Ocal_X \rightarrow \Ocal_W \rightarrow 0.\]
    In our local context, we are therefore interested in computing 
    \[H^{-q}(K^\bullet(\mathcal{G}\oplus \mathcal{F}, t \oplus 0)\otimes \mathcal{I}_W) \cong \bigoplus_{r+s = q}(H^{-r}(K^\bullet(\mathcal{G},t)\otimes \mathcal{I}_W)\otimes \bigwedge^s\mathcal{F}).\]
    Since $K^\bullet(\mathcal{G},t)$ is a resolution for $\Ocal_W$, we have that $H^{-r}(K^\bullet(\mathcal{G},t)\otimes \mathcal{I}_W) = Tor_r(\Ocal_W,\mathcal{I}_W)$. If we compute these by instead resolving the $\mathcal{I}_W$ by a truncated version of $K^\bullet(\mathcal{G},t)$, we see that $Tor_r(\Ocal_W,\mathcal{I}_W) \cong \bigwedge^{r+1}\mathcal{G}$. However, from the identification 
    \[\bigwedge^q(\mathcal{G}\oplus \mathcal{F}) \cong \bigoplus_{r+s = q}\bigwedge^r\mathcal{G}\otimes \bigwedge^s\mathcal{F}\]
    we see that
    \[H^{-q}(K^\bullet(\mathcal{G}\oplus \mathcal{F},t \oplus 0)\otimes \mathcal{I}_W) \cong \bigwedge^{q+1}(\mathcal{G}\oplus \mathcal{F})|_W/\bigwedge^{q+1}(\mathcal{F})|_W.\]
    This implies that the maps $H^{-q}\eta_w$ are injections and make the squares 
    \[
    \begin{tikzcd}
        &H^{-q}(\Delta^*(i_1\times i_2)_*\Ocal_{Y_1\times Y_2})\ar[d,"\cong"]\ar[r,"H^{-q}\eta_w"]&H^{-q}(w_*w^*\Delta_*(i_1\times i_2)_*\Ocal_{Y_1\times Y_2}) \ar[d,"\cong"]\\
        &w_*\bigwedge^q\mathcal{F}|_W \ar[r] &w_*\bigwedge^q(\mathcal{F}\oplus \mathcal{G})|_W
    \end{tikzcd}
    \]
    commute. Linking our two commuting squares together we have a commuting square 
    \[
    \begin{tikzcd}
        &H^{-q}(\Delta^*(i_1\times i_2)_*\Ocal_{Y_1\times Y_2})\ar[d,"\cong"]\ar[r,"H^{-q}\eta_w"]&H^{-q}(w_*w^*\Delta_*(i_1\times i_2)_*\Ocal_{Y_1\times Y_2}) \ar[d,"\varphi_q"]\\
        &w_*\bigwedge^q\mathcal{E}|_W \ar[r] &w_*\bigwedge^q(\mathcal{C}_{Y_1/X}\oplus \mathcal{C}_{Y_2/X})|_W
    \end{tikzcd}
    \]
    implying that our isomorphism is global and induced by $\varphi_q$.
\end{proof}

\section{Further directions}

Ideally, we would like to be able to make use of the technique of the proof of Theorem \ref{Theorem: Excess Intersection Formula} in cases where the intersection of the $Y_i$ is no longer an lci. However, following the logic of the proof, one would still need to be able to compare the local Koszul cohomologies to the restrictions of some global object in the case one is studying. Suppose that one found a global sheaf $\mathcal{G}$ on $X$ and morphism $f: \mathcal{G} \rightarrow w_*\bigwedge^q\bigoplus \mathcal{C}_{Y_i/X}|W$ such that there is an open affine cover of $X$ for which there are commuting diagrams on each element $U$ of the cover; 
\[
\begin{tikzcd}
        &H^{-q}(\Delta^*(i_1\times i_2)_*\Ocal_{Y_1\times Y_2})|_U\ar[d,"\cong"]\ar[r,"H^{-q}\eta_w|_U"]&H^{-q}(w_*w^*\Delta_*(i_1\times i_2)_*\Ocal_{Y_1\times Y_2}) |_U\ar[d,"\varphi_q|_U"]\\
        &\mathcal{G}|_U\ar[r, "f|_U"] &w_*\bigwedge^q(\mathcal{C}_{Y_1/X}\oplus \mathcal{C}_{Y_2/X})|_W.
    \end{tikzcd}
    \]
We want to show that the local vertical isomorphisms on the left glue up to give a global isomorphism $H^{-q}(\Delta^*(i_1\times i_2)_*\Ocal_{Y_1\times Y_2}) \cong \mathcal{G}$. The classical approach is to show that the restriction of the isomorphisms on two elements $U$ and $V$ of the cover agree on $U\cap V$. However, let $P$ be the categorical pullback in the diagram
\[
\begin{tikzcd}
        &P\ar[d,"\cong"]\ar[r]&H^{-q}(w_*w^*\Delta_*(i_1\times i_2)_*\Ocal_{Y_1\times Y_2}) \ar[d,"\varphi_q"]\\
        &\mathcal{G}\ar[r, "f"] &w_*\bigwedge^q(\mathcal{C}_{Y_1/X}\oplus \mathcal{C}_{Y_2/X})|_W.
    \end{tikzcd}
    \]
Since taking pullbacks commutes with restricting to open sets, by the universal property of pullbacks on $U\cap V$ we have two isomorphisms
\[H^{-q}(\Delta^*(i_1\times i_2)_*\Ocal_{Y_1\times Y_2})|_{U\cap V} \rightarrow P|_{U\cap V} \leftarrow H^{-q}(\Delta^*(i_1\times i_2)_*\Ocal_{Y_1\times Y_2})|_{U\cap V} \]
coming from the restrictions of the induced isomorphisms on $U$ and $V$. These isomorphisms have the property that when post-composed with the morphism
\[P_{U\cap V} \rightarrow H^{-q}(w_*w^*\Delta^*(i_1\times i_2)_*\Ocal_{Y_1\times Y_2})|_{U\cap V}\]
they are equal. To conclude that our isomorphisms agree then, it would be sufficient to know that our post-composition morphism is injective. This is what occurs in the result above. In cases where our post-composition morphism is not injective, the analysis of whether or not the local isomorphisms agree will have to be investigated more closely. 

\nocite{*}
\printbibliography
\end{document}